\documentclass[11pt,twoside,leqno]{aomamlt2e} 
  \pageno{1077}
\received{January 20, 2005}
 
\numberwithin{equation}{section}

\newtheorem{thm}{Theorem}%[section]

\newtheorem{lem}[thm]{Lemma}
\newtheorem{cor}[thm]{Corollary}

\newtheorem{prop}[thm]{Proposition}

\newtheorem{conj}[thm]{Conjecture}
\theorembodyfont{\upshape}

\newtheorem{defn}[thm]{{\it Definition}}

\newtheorem{say}[thm]{}
\newtheorem{exmp}[thm]{{\it Example}}

\newtheorem{ques}[thm]{{\it Question}}    %!!!!!!!!!!!!!!!!!!!!

\newtheorem{rem}[thm]{{\it Remark}}          

\newtheorem{defn-thm}[thm]{Definition-Theorem}  %!!!!!!!!!!!!!!!!!!!!!!!!
\newtheorem{defn-lem}[thm]{Definition-Lemma}  %!!!!!!!!!!!!!!!!!!!!!!!!
\theoremstyle{remark}
\newtheorem{claim}[thm]{{\it Claim}}

\renewcommand{\c}[0]{{\mathbb C}}  

\renewcommand{\o}[0]{{\mathcal O}} 
\newcommand{\z}[0]{{\mathbb Z}}

\renewcommand{\a}[0]{{\mathbb A}}

\newcommand{\p}[0]{{\mathbb P}}

\newcommand{\q}[0]{{\mathbb Q}}
\newcommand{\map}[0]{\dasharrow}
\newcommand{\qtq}[1]{\quad\mbox{#1}\quad}

\newcommand{\pic}[0]{\operatorname{Pic}}

\newcommand{\im}[0]{\operatorname{im}}    
\newcommand{\proj}[0]{\operatorname{Proj}}

\newcommand{\ext}[0]{\operatorname{Ext}}    
\newcommand{\Hom}[0]{\operatorname{Hom}}

\newcommand{\aut}[0]{\operatorname{Aut}}

\newcommand{\sing}[0]{\operatorname{Sing}}    
\newcommand{\ex}[0]{\operatorname{Ex}}

\newcommand{\hilb}[0]{\operatorname{Hilb}}

\newcommand{\ns}[0]{\operatorname{NS}}
\newcommand{\PGL}[0]{\operatorname{PGL}}
\newcommand{\GL}[0]{\operatorname{GL}}
\newcommand{\SO}[0]{\operatorname{SO}}
\newcommand{\Ogp}[0]{\operatorname{O}}

\def\into{\DOTSB\lhook\joinrel\to}

\begin{document}
\currannalsline{164}{2006} 

 \title{Non-quasi-projective moduli spaces}
\author{J\'anos Koll\'ar}

 \institution{Princeton University, Princeton, NJ\\
\email{kollar@math.princeton.edu}}

\centerline{\bf Abstract}\vglue8pt

We show that every smooth toric variety                              
(and many  other algebraic spaces as well)  can be realized                     
as a moduli space for smooth, projective, polarized varieties.                  
Some of these are not quasi-projective.                                        
This contradicts a recent paper                                                 
(Quasi-projectivity of moduli spaces of polarized varieties,  {\it Ann.\ of
Math\/}.\     
 {\bf 159}  (2004) 597--639.).   

\vglue8pt

A polarized variety is a pair $(X,H)$ consisting of a smooth
projective variety $X$ and a linear equivalence class of
 ample divisors $H$ on $X$.
For simplicity, we look at the case when $X$ is smooth, 
numerical and linear equivalence coincide for divisors on $X$, 
$H$ is very ample
and $H^i(X,\o_X(mH))=0$ for $i,m>0$.
A well established route to construct moduli spaces
of such pairs is to embed $X$ into $\p^N$ by $|H|$.
The pair $(X,H)$ and the embedding $X\into \p^N$
determine each other up to the action of $\PGL(N+1)$.
Deformations of $(X,H)$  cover an open subset $U(X,H)$ of
the Hilbert scheme $\hilb(\p^N)$ with Hilbert polynomial $\chi(X,\o_X(mH))$.
One can then view the quotient $U(X,H)/\PGL(N+1)$ as the moduli space
of the pairs $(X,H)$.
(See \cite[App.~5]{git} or \cite[Ch.~1]{vie} for general
introductions to moduli problems.)

The action of $\PGL(N+1)$ can be  bad along some orbits,
and therefore one has to make  additional assumptions to ensure
that the quotient $U(X,H)/\PGL(N+1)$ is reasonable.
 The optimal condition seems to be
to require that the action be {\it proper}. This is
equivalent to assuming that $U(X,H)/\PGL(N+1)$ exists as a
separated complex space or as a separated algebraic space
 \cite{koll-quot}, \cite{ke-mo}.

A difficult result of Viehweg (cf.\ \cite{vie}) shows that
if the canonical class $K_X$ is assumed nef then
$U(X,H)/\PGL(N+1)$ is a quasi-projective scheme.

A recent paper \cite{sch-tsu} asserts the quasi-projectivity 
of moduli spaces of
polarized varieties for  arbitrary
$K_X$, whenever the quotient $U(X,H)/\PGL(N+1)$ exists as a
separated
algebraic space.

The aim of the present note is to confute this claim.
The examples (\ref{toric.as.mod.cor}) and  (\ref{versal.mod.exmps}) show
that the quotients $U(X,H)/\PGL(N+1)$ can contain smooth,
proper subschemes which are not projective.

In the examples $X$ is always a rational variety,
but there are many more such cases as long as
$X$ is ruled. This leaves open the question of
quasi-projectivity of the quotients
$U(X,H)/\PGL(N+1)$ when $X$ is not uniruled but $K_X$ is not nef.

We work over an algebraically
closed field of characteristic zero,
though some of the examples apply in any characteristic.
\vglue-20pt
 \phantom{up}
\section{First examples}

\vglue-18pt
 \phantom{up}
\begin{say}[Versions of quasi-projectivity for moduli functors]\label{vers.qp}
In asserting
that certain moduli spaces are quasi-projective,
one hopes to show that
an algebraic space $S$ is quasi-projective
if $S$ ``corresponds'' to a family of pairs $(X,H)$
in our class. 
There are at least three ways to formulate 
a precise meaning of ``corresponds''.
(In order to avoid  scheme theoretic complications,
let us assume that $S$ is normal.)
 
\vskip2pt
(\ref{vers.qp}.1. There is a family over $S$.) That is, 
there is   a smooth, proper morphism
of algebraic spaces $f:U\to S$  and an
$f$-ample Cartier divisor $H$ such that
 every fiber $(U_s,H|_{U_s})$ is in our class
and $(U_s,H|_{U_s})\cong (U_{s'},H|_{U_{s'}})$ if and only if $s=s'$.

\vskip2pt
(\ref{vers.qp}.2. There is a family over some scheme over $S$.) That is, 
there are  a surjective and open morphism $h:T\to S$,  
 a smooth, proper morphism
of algebraic spaces $f:U\to T$  and an
$f$-ample Cartier divisor $H$ such that
 every fiber $(U_t,H|_{U_t})$ is in our class
and $(U_t,H|_{U_t})\cong (U_{t'},H|_{U_{t'}})$ if and only if $h(t)=h(t')$.
(One can always reduce to the case when $h:T\to S$
is the geometric quotient by a $\PGL$-action, but
in many constructions quotients by smaller groups appear naturally.)

\vskip2pt

(\ref{vers.qp}.3. There is a universal family over some scheme over $S$.)
That is, we have  $h:T\to S$  and $f:U\to T$ 
as in (\ref{vers.qp}.2) but we also assume that every local deformation of a
polarized fiber  $(U_t,H|_{U_t})$ is induced
from $f:U\to T$. 
\end{say} 
\phantom{up}
\vglue-22pt

 All the 
approaches to quasi-projectivity of quotients that I know of
work equally well for any of the three cases.
(For instance, although the main assertion \cite[Thm.~1]{sch-tsu} 
explicity assumes
local versality as in (\ref{vers.qp}.3),
 the key technical steps
\cite[Thms.~4, 5]{sch-tsu} assume only the more general setting
of (\ref{vers.qp}.2).)
Nonetheless, a counterexample to the variant (\ref{vers.qp}.2)
need not yield automatically a 
counterexample in the setting of (\ref{vers.qp}.3).

I  start with  examples as  in  (\ref{vers.qp}.2) where quasi-projectivity
fails; these are  the
weak examples (\ref{weak.exmp}).
Then we analyze  deformations of some of these polarized pairs to 
show that  quasi-projectivity also
fails under the assumptions  of  (\ref{vers.qp}.3).
These examples are given in Section 4.

\vglue-22pt \phantom{up}
\begin{say}[Weak examples]\label{weak.exmp} 
Let $W^0$ be a smooth, quasi-projective variety
of dimension at least 2
and $G$ a reductive algebraic group acting  on $W^0$.
Let $W\supset W^0$ be a $G$-equivariant compactification of $W^0$.

The moduli space of the pairs $(w,W)$ consisting of $W$ 
(thinking of it as fixed)
and a variable point $w\in W^0$ is naturally a quotient of 
$W^0/G$. These pairs can also be identified with
pairs  $(B_wW,E)$ where $E\subset B_wW$ is the exceptional divisor
of the blow up $\pi_w:B_wW\to W$ of $w\in W$. 
 Fix a sufficiently ample 
$G$-invariant linear equivalence class of divisors $H$ on $W$.
Then $H_w=\pi^*H-E$ is ample on $B_wW$ and 
$(B_wW,H_w)$ uniquely determines $(B_wW,E)$ (cf.\ (\ref{rigidify.X}).

Thus we  obtain  a $G$-equivariant morphism of   $W^0$
to the moduli space of the
polarized pairs $(B_wW,H_w)$.

Assume now that in the above example the
following conditions are satisfied: 
\vskip4pt

 (\ref{weak.exmp}.1) the $G$-action is proper on  $W^0$, 

\vskip4pt (\ref{weak.exmp}.2) $W^0/G$ is not quasi-projective, and

\vskip4pt (\ref{weak.exmp}.3) $\aut(W)=G$.
\vskip4pt

The quotient  $W^0/G$ exists as an algebraic space
by the general quotient results of  \cite{koll-quot}, \cite{ke-mo}. 
In (\ref{vers.qp}.2) set $S=W^0/G$ and $T=W^0$.
The pairs $(B_wW,H_w)$ give a family of
polarized varieties over $W^0$.
Furthermore, isomorphisms between
two polarized pairs $(B_wW,H_w)$ and $(B_{w'}W,H_{w'})$
correspond to  isomorphisms between
 the pairs $(w,W)$ and $(w',W)$, 
and by (\ref{weak.exmp}.3), these in turn are given by those
elements of $G$ that map $w$ to $w'$.
In particular, two polarized pairs $(B_wW,H_w)$ and $(B_{w'}W,H_{w'})$
are isomorphic if and only if $w$ and $w'$ are in the same $G$-orbit.

Thus we have realized the non-quasi-projective
algebraic space  $W^0/G$ as a moduli
space of smooth, polarized varieties in the sense of (\ref{vers.qp}.2).
\end{say}
\phantom{up}
\vglue-20pt
Now we must find examples where the three conditions
of (\ref{weak.exmp}) are satisfied.
We start by reviewing some of the known examples
of proper $G$-actions with non-quasi-projective
quotient. The condition $\aut(W)=G$ should hold
for most $G$-equivariant compactifications, but
it will take some effort to prove that such a $W$
exists in many cases.

\vglue-20pt
\phantom{up}
\begin{say}[Examples of non-quasi-projective quotients]
\label{nonqp.quot.exmp}
There are many examples of  $G=\PGL$ or  a torus $G=(\c^*)^m$ acting 
properly on a
smooth quasi-projective variety $W^0$ such that
$W^0/G$ is not quasi-projective.

Here we show two examples where a torus or $\PGL(n)$ acts 
properly on
an open subset of projective space and the
quotient is smooth, proper but not projective in the torus case
and a smooth algebraic space which is  not a scheme in the $\PGL(n)$ case.

\vglue3pt
(\ref{nonqp.quot.exmp}.1) By a result of \cite[Thm.~2.1]{cox}, every smooth 
toric variety
can be written as the geometric quotient of an open subset
$U\subset \c^N\setminus\{0\}$ by a suitable subtorus of
$(\c^*)^N$. There are many proper but
nonprojective toric varieties (see, for instance, \cite[\S 2.3]{oda}),
and so we have our first set of examples.

(\ref{nonqp.quot.exmp}.2) Here we work with $\PGL(3)$, but the
construction can be generalized to any $\PGL(n)$ for $n\geq 3$.

Fix  $d$ and let
$U_d\subset |\o_{\p^2}(d)|$ be the open set consisting of curves $C$
such that 
\begin{enumerate}
\item[(i)] $C$ is smooth, irreducible and 
the genus of its normalization is 
$>\frac12\binom{d-1}{2}$.
\item[(ii)] $C$ is not fixed by any of the automorphisms of $\p^2$.
\end{enumerate}
We claim that $\aut(\p^2)$ operates properly and  freely
on $U_d$. Indeed, the action is set theoretically free by (ii).
Properness is equivalent to uniquenes of specialization:

\begin{claim} Let $S$ be the spectrum of a DVR.
A family of smooth plane curves of degree $d$ 
over the generic point 
$S^*\subset S$ has at most one extension to a family over $S$
where the central fiber is in $U_d$.
\end{claim}

\Proof Assume that we have  a family $X^*\to S^*$
and two extensions  $X_1,X_2\to S$
with central fibers $C_1,C_2$. 
If the natural map $X_1\map X_2$ is an isomorphism at the
generic point of $C_1$, then the two families are isomorphic
by (\ref{mats-mum.lem}).

Otherwise, let
$Y\to S$ be the normalization of the main component of the
fiber product $X_1\times_SX_2$.
The central fiber of $Y\to S$ dominates both $C_1,C_2$,
hence it has two irreducible components, both of geometric genus
$>\frac12\binom{d-1}{2}$. Thus the
sum of the geometric genera of the irreducible components
of the central fiber is bigger than the
geometric genus of the generic fiber, a contradiction.\Endproof\vskip4pt

Let us consider a general curve $C\subset \p^2$
which has multiplicity $\geq m$ at a given point $p\in \p^2$.
Our condition for the  geometric genus is
$$
\binom{d-1}{2}-\binom{m}{2}>\frac12\binom{d-1}{2},
$$
which is asymptotically equivalent to $m<d/\sqrt{2}$.

On the other hand, if $m>2d/3$ and $p=(0:0:1)$ then the subgroup
$(t,t,t^{-2})$ shows that $[C]\in |\o_{\p^2}(d)|$
is unstable. Since $2/3<1/\sqrt{2}$, we obtain:

\begin{claim} For large $d$, there are curves $C$ with 
 $[C]\in U_d$ such that $[C]$ is unstable.\hfill\qed
\end{claim}

\begin{cor} For large $d${\rm ,} the quotient 
$U_d/\aut(\p^2)$ is a smooth algebraic space which is not a scheme.
\end{cor}

\Proof The quotient 
$U_d/\aut(\p^2)$ is a smooth algebraic space by 
 \cite{koll-quot}, \cite{ke-mo}. 
Let $\pi:U_d\to U_d/\aut(\p^2)$ denote the quotient map.

Pick a curve 
 $[C]\in U_d$ such that $[C]$ is unstable.
We claim that $[C]\in U_d/\aut(\p^2)$
has no neighborhood which is affine.
Indeed, if $W\subset U_d/\aut(\p^2)$
is any quasi-projective subset, then by \cite[Converse 1.12]{git},
 its preimage
$\pi^{-1}(W)\subset |\o_{\p^2}(d)|$
consists of semi-stable points 
with respect to some polarization on $|\o_{\p^2}(d)|\cong \p^N$.
Since $|\o_{\p^2}(d)|$ is a projective space
and $\aut(\p^2)=\PGL(3)$ has no nontrivial homomorphisms to $\c^*$,
up to powers one has only the standard polarization, and so
$\pi^{-1}(W)\subset |\o_{\p^2}(d)|$
consists of semi-stable points 
with respect to the usual polarization.
Thus $\pi^{-1}(W)$ cannot contain
$[C]$ since $C$ is unstable.\hfill\qed
\end{say}

The third  requirement  (\ref{weak.exmp}.3)
is to find a compactification of a $G$-variety 
whose automorphism group is exactly $G$.
Thus we need to consider the following general problem.

\begin{ques}\label{rigid.comp.quest} Let $G$ be an algebraic group acting
on a quasi-projective variety $W^0$. When can one find
a projective compactification  $W^0\subset W$ such that
$\aut(W)=G$?
\end{ques}

There are some cases when this cannot be done.
The simplest counterexample occurs when $W^0$ is
projective; here we have no choices for $W$.
The answer can be negative even if $W^0$ is affine.
For instance, consider the action of
$\Ogp(n)$ on $W=\p^{n-1}$. Here there are only two orbits;
let $W^0$ be the open one. As the complement $W\setminus W^0$
is a single codimension 1 orbit, there are no 
$\Ogp(n)$-equivariant blow ups to make, so
$W=\p^{n-1}$ is the unique $\Ogp(n)$-equivariant compactification
of $W^0$ and $\aut(W)=\PGL(n)$ is bigger than $\Ogp(n)$.

The question becomes more reasonable if we assume that
$G$ acts properly on $W^0$. There is still an easy negative example,
 $G=W^0=\c^*$, but there may not be any others where $G$ is reductive.
In the next two sections, we prove the following partial result.
 \newcommand{\ee}{{\hskip1pt\rm \'{\hskip-6.5pt \it e}}}

\begin{prop} \label{rgid.comp.prop}
 Let $G$ be either a torus $(\c^*)^n$ or $\PGL(n)$.
Let $W^0$ be a smooth variety with a generically free and  proper $G$-action
such that  $\rho(W^0)\!=\!0${\rm ,} that is{\rm ,} $\pic(W^0)$ is a torsion group.
 Assume that there is a
\/{\rm (}\/not necessarily $G$-equivariant\/{\rm )}\/ smooth compactification
$W^0\subset W^*$ such that its
N\ee ron-Severi group $NS(W^*)$ is $\z$.

Then there is a smooth $G$-equivariant compactification
$W\supset W^0$ and an ample divisor class $H$
such that $\aut(W,H)=G$.

Moreover{\rm ,} if $W'\to W$ is any other $G$-equivariant compactification
dominating $W$ then there is an ample divisor class $H'$
such that $\aut(W',H')=G$.
\end{prop}

Putting this together with
(\ref{nonqp.quot.exmp}.1) we obtain the following:

\begin{cor}\label{toric.as.mod.cor} Every smooth toric variety
can be written as a moduli space of
smooth{\rm ,} polarized varieties as in {\rm (\ref{vers.qp}.2).}\hfill\qed
\end{cor}

By a theorem of \cite{wlo}, a smooth proper variety $X$
can be embedded into a smooth toric variety if and only if every two points of $X$ are contained in an open affine subset.
Thus (\ref{toric.as.mod.cor}) implies that 
a smooth proper variety $X$
can be written as a moduli space of
smooth, polarized varieties as in (\ref{vers.qp}.2)
provided every two points of $X$ are contained in an open affine subset.

In the next section we start the proof of of (\ref{rgid.comp.prop})
by  finding $W$ such that the connected component of
$\aut(W)$ is $G$. After that we choose the polarization
$H$ such that $\aut(W,H)$ equals the connected component of
$\aut(W)$.

\section{Rigidifying by compactification}

\begin{defn} Let $X$ be a proper variety and $\ns(X)$ its
N\'eron-Severi group. The automorphism group
$\aut(X)$ acts on $\ns(X)/(\mbox{torsion})$; let $\aut^0(X)$ denote the
kernel of this action.
\end{defn}

\begin{lem}\label{aut0.in.birmaps}
 Let $f:Y\to X$ be a proper{\rm ,} birational morphism between
smooth projective varieties. Then $\aut^0(Y)\subset \aut^0(X)$.
\end{lem}

\Proof The exceptional set $\ex(f)$ is a union of divisors
and an exceptional divisor is not linearly equivalent to any other
effective divisor. Thus 
$\aut^0(Y)$ stabilizes $\ex(f)$ and so every $\sigma\in \aut^0(Y)$
descends to an automorphism $\sigma_X$ of $X\setminus f(\ex(f))$.
Since $f(\ex(f))$ has codimension at least 2 and 
 $\sigma_X$  fixes an ample divisor,  $\sigma_X\in \aut(X)$
by (\ref{mats-mum.lem}). \hfill\qed

\begin{lem}[\cite{ma-mu}]\label{mats-mum.lem}
 Let $X,X'$ be normal{\rm ,} projective varieties
and\break $Z \subset X${\rm ,} $Z'\subset X'$  closed subsets of codimension $\geq 2$.
Let $\phi:X\setminus Z\to X'\setminus Z'$
be an isomorphism. Assume that there are ample divisors
$H$ on $X$ and  $H'$ on
$X'$ such that $\phi^{-1}(H')=H$.
Then $\phi$ extends to an isomorphism
$\Phi:X\to X'$.\hfill\qed
\end{lem}

We deal with the difference between $\aut^0(X)$ and $\aut(X)$ later.
Now we concentrate on answering (\ref{rigid.comp.quest})
for certain cases that are of special interest in moduli
constructions. To this end we introduce another subgroup of $\aut$.

\begin{defn}  Let $W^0$ be a variety with a $G$-action.
and  $W\supset W^0$ a  $G$-equivariant compactification.
Let $\aut_{\partial}(W)\subset \aut(W)$ be the subgroup consisting
of all automorphisms which stabilize every $G$ orbit in $W\setminus W^0$.
\end{defn}

\begin{lem} \label{aut-partial.in.aut0}
 Let $W^0$ be a variety with a $G$-action{\rm ,} $G$ connected.
Let  $W\supset W^0$ be a  $G$-equivariant smooth compactification.
If $\rho(W^0)=0$ then $\aut_{\partial}(W)\subset \aut^0(W)$.
\end{lem}

\Proof Since  $\rho(W^0)\!=\!0$, the divisorial irreducible components
of $W\backslash W^0$ generate $\ns(W)_{\q}$. Since $G$ is connected,
each irreducible component of\break $W\backslash  W^0$ is fixed by $G$, hence by
$\aut_{\partial}(W)$. Thus $\aut_{\partial}(W)$ acts trivially
on $\ns(X)/(\mbox{torsion})$.\hfill\qed

\begin{cor} \label{aut-partial.in.birmaps}
Let $W^0$ be a variety with a $G$-action{\rm ,} $G$ connected.
Let   $W_i\supset W^0$ be   $G$-equivariant smooth compactifications
and $W_1\to W_2$ a proper{\rm ,} birational  $G$-equivariant morphism.
If $\rho(W^0)=0$ then $\aut_{\partial}(W_1)\subset \aut_{\partial}(W_2)$.
\end{cor}

\Proof From (\ref{aut-partial.in.aut0})
we know that $\aut_{\partial}(W_1)\subset \aut^0(W_1)$
and $\aut^0(W_1)\subset \aut^0(W_2)$ by
(\ref{aut0.in.birmaps}). Since every $G$-orbit in
$W_2$ is the image of a $G$-orbit in $W_1$, 
the inclusion $\aut_{\partial}(W_1)\subset \aut_{\partial}(W_2)$
follows.\hfill\qed

\begin{exmp} It is worth noting that (\ref{aut-partial.in.birmaps})
can fail if $\rho(W^0)>0$. Start with the $\Ogp(4)$ action on
$W^0=(xy-uv=0)\setminus\{(0,0,0,0)\}\subset \a^4$.
Let $W\subset W^0$ be its closure in $\p^4$.
Let $W_1\to W$ be the blow up of the origin and $W_2\to W$
the blow up of $(x=u=0)$. The induced map
$W_1\to W_2$ is a blow up of a single smooth rational curve.
$\Ogp(4)$ acts on $W_1$ but only $\SO(4)$ acts by automorphisms on $W_2$.
The involution $(x,y,u,v)\mapsto (x,y,v,u)$
lifts to a birational involution on $W_2$ which is not
an automorphism.
\end{exmp}

\begin{prop}\label{aut-partial=aut0-stably}
 Let $G$ be a connected algebraic group and 
 $W^0$  a smooth variety with a  $G$-action
such that  $\rho(W^0)=0$ and $\dim G\leq \dim W^0-2$. 
Then there is a smooth $G$-equivariant compactification
$W\supset W^0$ such that $\aut_{\partial}(W)=\aut^0(W)$.

Moreover{\rm ,} if $W'\to W$ is any other $G$-equivariant compactification
dominating $W$ then $\aut_{\partial}(W')=\aut^0(W')$.
\end{prop}

\Proof Let us start with any smooth 
$G$-equivariant compactification $W_1\supset W^0$.
As  $\aut^0$ can only decrease
under further blow ups, we can assume that it is already  minimal.
That is, if 
$W'\to W$ is any other $G$-equivariant compactification
then  $\aut^0(W')=\aut^0(W)$.

Assume now that $\aut_{\partial}(W)\neq \aut^0(W)$. Then there are a
$\sigma\in \aut^0(W)$ and a $G$-orbit $Z\subset W\setminus W^0$ 
such that $\sigma(Z)\neq Z$. After some preliminary\break $G$-blow ups we can
 blow up $Z$ to get $W_Z\to W$. 
Since $\dim G\leq \dim W^0-2$, this blow up is nontrivial and
the preimage of $Z$ is an exceptional divisor $E_Z$.
We also know that $E_Z$ is not numerically equivalent to
any other effective divisor and 
it is not stabilized by $\sigma$. Thus
$\aut^0(W_Z)\neq\aut^0(W)$, a contradiction. \hfill\qed

\begin{say}[First examples with $G=\aut_{\partial}(W)$]\label{2exmps}
\vskip4pt

(\ref{2exmps}.1) Let $G=(\c^*)^n$ be the torus with its left action on itself.
A natural compactification is $W=\p^n$. The coordinate ``vertices''
are fixed by $G$ and by no other automorphism of $W$.
Thus $G=\aut_{\partial}(W)$.
Moreover, if 
 $W'\to W$ is any other $G$-equivariant compactification
dominating $W$ then $G\subset \aut_{\partial}(W')\subset \aut_{\partial}(W)$;
hence $G= \aut_{\partial}(W')$.
\vskip4pt

(\ref{2exmps}.2) Let $G=\PGL(n)$  with its left action on itself.
A natural compactification is $W=\p(M_n)$ coming from the $\GL(n)$ action
on $n\times n$-matrices by left multiplication.
The $(n-1)$-dimensional $G$-orbits are of the form
 $\p^{n-1}\times (a_1,\dots,a_n)$ where we think of the points in $\p^{n-1}$ 
as column vectors. The union of all 
$(n-1)$-dimensional $G$-orbits is $\p^{n-1}\times \p^{n-1}$ under the Segre
 embedding. From this we conclude that $\aut_{\partial}(W)$
acts on $\p^{n-1}\times \p^{n-1}$ as multiplication on the first factor.
Since the image of  $\p^{n-1}\times \p^{n-1}$ under the Segre
 embedding is not contained in any hyperplane, this implies that
$\aut_{\partial}(W)=\PGL(n)$.

As before,  if 
 $W'\to W$ is any other $G$-equivariant compactification
dominating $W$ then  $\aut_{\partial}(W')=\PGL(n)$ as well.
\end{say}

We are now ready to to answer (\ref{rigid.comp.quest})
for $(\c^*)^n$ and for $\PGL(n)$.

\begin{prop} Let $G$ be either  $(\c^*)^n$ or $\PGL(n)$.
Let $W^0$ be a smooth variety with a
generically free  and  proper $G$-action
such that  $\rho(W^0)=0$ and $\dim G\leq \dim W^0-2$. 
Then there is a smooth $G$-equivariant compactification
$W\supset W^0$ such that $\aut^0(W)=G$.

Moreover{\rm ,} if $W'\to W$ is any other $G$-equivariant compactification
dominating $W$ then $\aut^0(W')=G$.
\end{prop}

\Proof Let $Z\supset W^0/G$ be any compactification
and choose any compactification $W_1\supset W^0$
such that there is a morphism $h:W_1\to Z$.
By further $G$-equivariant blow ups in $W_1\setminus W^0$,
(\ref{aut-partial=aut0-stably}) gives $W_2\supset W^0$ such that
$\aut_{\partial}(W_2)=\aut^0(W_2)$ and neither of these groups
changes under further $G$-equivariant blow ups in $W_2\setminus W^0$.

Pick a big linear system of Weil divisors $|B|$ on $Z$
and let $|M|$ be the moving part of the linear system given by
a pull back of the general member of $|B|$.
Then $|M|$ is a linear system which gives the map 
$h:W_2\to Z$
over some open subset of $Z$. ($Z$ is not projective in general,
and may not even have any Cartier divisors. That is why we have to
find $|M|$ in this roundabout way.)

Any element of $\aut^0(W_2)$ sends $|M|$ to itself,
hence $h:W_2\to Z$ is $\aut^0(W_2)$-equivariant.

General fibers of $h$ contain a  $G$-orbit which is in 
$W_2\setminus W^0$, and so every general fiber of 
$h$ is $\aut^0(W_2)$-stable since
$\aut_{\partial}(W_2)=\aut^0(W_2)$.

Pick any $\sigma\in \aut^0(W_2)$ and look at its action
$\sigma_z$ on $h^{-1}(z)$ for general $z\in Z$.

Since $z$ is general, the fiber $h^{-1}(z)$
is a smooth projective $G$-equivariant compactification
of $G$ acting on itself. We claim that
$\sigma_z=g(z,\sigma)$ for some $g(z,\sigma)\in G$.
This follows from (\ref{2exmps})
if $h^{-1}(z)$ dominates the  compactifications considered there.
Otherwise, by further blow ups we could get $W_3\to W_2$ such that
the birational transform of 
$h^{-1}(z)$ dominates 
the standard compactifications considered in  (\ref{2exmps}).
This would, however, mean that $\sigma$ does not lift to
$\aut^0(W_3)$, a contradiction to our assumption that
$\aut^0(W_2)$  does not
change under further $G$-equivariant blow ups.

Thus we conclude that $G$ and $G':=\aut^0(W_2)$
both act on $W_2$ in such a way that for a general
$w\in W^0$,
\begin{enumerate}
\item $Gw=G'w$, and
\item the $G'$-action on $Gw$ is via a homomorphism
$\rho_w:G'\to G$.
\end{enumerate}
Let $H'_w\subset G'$ be the kernel of $\rho_w$.
Since $G$ is reductive, $H'_w$ contains the
unipotent radical  $U'\subset G'$.
The quotients $H'_w/U'$ are normal subgroups 
of the reductive group $G'/U'$, and they depend continuously
on $w$ over an open set of $W$ (\ref{stabs.cont.dep}).
A continuously varying family of normal subgroups would give a
continuously varying family of finite dimensional
representations, but a reductive group has
only discrete series representations
in finite dimensions.
This implies that $H'_w$ is independent of $w$ for general
$w\in W$ and so the\break $H'_w$-action is trivial on $W_2$.
But $H'_w\subset \aut^0(W_2)$, thus $H'_w$ is the trivial group
and so $G=\aut^0(W_2)$.\hfill\qed

\begin{exmp} The example $G'=\c^2_{x,y}$ acting on
$\c^2_{u,v}$ as $$(u,v)\mapsto (u,v+x-uy)$$ shows that
the above argument does not work if $G$ is not reductive.
\end{exmp}

\begin{rem} \label{stabs.cont.dep} Let $G$ be an algebraic group
acting on a variety $X$. The stabilizer subgroups $G_x$ of points $x\in X$
are the same as the fibers of $G\times X\to X\times X$ over the diagonal.
Thus we see that
\begin{enumerate}
\item  the dimension of $G_x$ is a constructible function
on $X$,
\item  the number of connected components of $G_x$ is a constructible function
on $X$,
\item  the subgroups $G_x\subset G$ depend continuously on $x$
for $x$ in a suitable open subset of $X$.
\end{enumerate}
\end{rem}

\section{Rigidifying using polarizations}

Let $X$ be a proper variety and $H$ an ample divisor
on $X$. Then $\aut(X,H)$ can be viewed as  a closed subgroup
of $\p(H^0(X,\o_X(mH)))$ for $m\gg 1$. Hence
$\aut(X,H)$ is an algebraic group and so it has only finitely
many connected components. This implies that the action
of $\aut(X,H)$ on $\ns(X)$ is through a finite group.

While not crucial, it will be convenient for us to choose a
polarization such that $\aut(X,H)$ acts trivially on
$\ns(X)$. In particular,  $\aut(X,H)=\aut^0(X)$.

\begin{lem}\label{rigidify.X}
 Let $g:Y\to X$ be a birational morphism between smooth{\rm ,}
projective varieties. Assume that  $\ns(X)\cong \z$.
Then there is an ample divisor $H^*$ on $Y$ such that
$\aut(Y,H^*)=\aut^0(Y)$.
\end{lem}

\Proof Let $g:Y\to X$ be a birational morphism between smooth,
projective varieties with exceptional divisors
$E_i$. Let $H$ be any ample divisor on~$X$.

Let $H_Y$ be ample on $Y$. Then we can write
$H_Y=g^*(g_*H_Y)-\sum a_iE_i$ for some $a_i>0$.
For reasons that will become clear soon, let us change the
$a_i$ a little so that
we get an ample $\q$-divisor
$H'_Y:=g^*(g_*H_Y)-\sum a'_iE_i$ where  the $a'_i$ are different from
each other.
Choose $m$ such that $mH-g_*H'_Y$ is ample
on~$X$. Then 
$$
mg^*H-\sum a'_iE_i=g^*(mH-g_*H'_Y)+H'_Y\qtq{is ample  on $Y$.}
$$
Let us multiply through with the common denominator of the $a'_i$
to get natural numbers $b_i$ and $m_0$ such that
$$
H_m:=mg^*H-\sum b_iE_i\qtq{is ample for $m\geq m_0$,}
$$
and the $b_i$ are different from each other.

Write $K_Y=f^*K_X+\sum e_iE_i$, where $e_i>0$ for every $i$.
Choose a natural number $c$ such that $ce_i-b_i\geq 0$ for every $i$.
Finally, choose $m_1$ such that $mH+cK_X$ is very ample
on $X$ for $m\geq m_1$.

\begin{claim} For $m\geq \max\{m_0,m_1\}$, the polarized variety $(Y,H_m)$
uniquely determines $f:Y\to X$ and also $\sum b_iE_i$.
\end{claim}

\Proof Given $H_m$, we consider the linear system
$$
\begin{array}{rcl}
|H_m+cK_Y|&=&|mg^*H-\sum b_iE_i+cg^*K_X+\sum ce_iE_i|\\
&=&
g^*|mH+cK_X|+\sum (ce_i-b_i)E_i,
\end{array}
$$
where the second equality holds 
since an effective exceptional divisor does not increase
a linear system that is pulled back from the base.
As $mH+cK_X$ is very ample by assumption, we see
that we recover $g:Y\to X$ as
given by $|H_m+cK_Y|$.
Furthermore, since the fixed part $\sum (ce_i-b_i)E_i$
 is also determined by  $|H_m+cK_Y|$,
we also recover  $\sum b_iE_i$.\hfill\qed
\medskip

Now we use the fact that all the $b_i$ are different from each other.
This implies that every automorphism of $(Y,H_m)$
maps each $E_i$ to itself. Furthermore, $g^*H$ is
also mapped to itself. Since $X$ has Picard number 1,
these  together generate  a finite index subgroup of 
the free abelian group $\ns(Y)$.
Thus $\aut(Y,H_m)$ acts trivially on $\ns(Y)$.\hfill\qed

\section{Locally versal examples}

Start with  $\a^{n}$ with the standard $(\c^*)^{n}$-action.
Let $T\subset (\c^*)^{n}$ be a subtorus and 
 $U\subset \a^{n}$  a $(\c^*)^{n}$-invariant
open set on which $T$ acts properly.
In (\ref{weak.exmp}) we showed how to construct a
moduli problem for smooth polarized varieties whose
moduli space is $U/T$. These give examples of
moduli spaces as in (\ref{vers.qp}.2), but
in general the local versality condition of
(\ref{vers.qp}.3) fails.

In this section we  present a version of the construction
where we can control local versality as well.
The key point is to get a rather explicit
series of examples as in 
(\ref{weak.exmp}) where we can describe all deformations
in a uniform way.

This should be possible to do
in most cases, but the combinatorial aspects
of finding explicit resolutions and describing their deformations seem
rather daunting.
So here I consider a class of special examples, where
the ancillary problems are easy to handle.

\begin{say}[Conditions on $T$]\label{torus.conds} We  assume from now on that
our torus $T$ and its action on affine space is of the following form.

Start with $\a^{s+t}$ with the standard $(\c^*)^{s+t}$-action.
Fix positive integers $c_{ij}$ and let 
$T=T(c_{ij})=\im[(\c^*)^{t}\to  (\c^*)^{s+t}]$ where the map is given by
$$
(\lambda_1,\dots,\lambda_t)
\mapsto \bigl(\textstyle{\prod_j}\lambda_j^{c_{1j}},\cdots,
\textstyle{\prod_j}\lambda_j^{c_{sj}},
\lambda_1,\dots,\lambda_t\bigr).
$$
Let $U\subset \a^{s+t}$ be a $(\c^*)^{s+t}$-invariant
open set on which $T$ acts properly.
\end{say}

\begin{say}[Choosing singular moduli problems] \label{csmp.say}

Take $\p^n\supset \a^{n}$ with coordinate vertices
$p_0={\mathbf 0}\in \a^n$ and $p_1,\dots,p_n$ at infinity.

Given $n=s+t$  as in (\ref{torus.conds}) and a positive integer $d$,
let $\mathcal X=\mathcal X(s,t,d)$ be the set of all 
varieties $X$  that are obtained as $f:X\to \p^n$ from $\p^n$ by performing
\begin{enumerate}
\item a weighted blow up (see (\ref{w.bu.defn}))
with weights $(d^s,1^t)$ at $p_0$,
\item ordinary blow ups 
at the $n$ points $p_1,\dots,p_n$ and at a further point $q\in U$.
\end{enumerate}
Let $E_0,\dots,E_n,E_{n+1}$ be the corresponding exceptional divisors.
As in (\ref{rigidify.X}) choose a polarization $H$ on $X$ of the form
$f^*\o_{\p^n}(m)-\sum b_iE_i$ such that
the map $f:X\to \p^n$ and the $E_i$ are determined by
the pair $(X,H)$.
We obtain the set of polarized pairs $(\mathcal X,\mathcal H)$

Weighted blow ups depend on the choice of a
local coordinate system, and for 
weights $(d^s,1^t)$ we show that they correspond uniquely
to certain ideals  $I_d\subset \o_{{\mathbf 0},\a^{s+t}}$.
Thus $\mathcal X$ has a natural scheme structure as
a subset of the Hilbert scheme of points
$\hilb(\a^{s+t})$ corresponding to the union of
$q\in U$ and $\o_{{\mathbf 0},\a^{s+t}}/I_d$.

\begin{prop}[Notation as above] Assume that $t\geq 3$. Then
 $(\mathcal X,\mathcal H)$ is locally versal
and the isomorphisn classes of the
polarized pairs $(X,H)\in (\mathcal X,\mathcal H)$
are in one-to-one correspondence with the
$(\c^*)^{s+t}$-orbits on $\mathcal X$.
\end{prop}

\Proof Every defomation of a smooth point blow up is 
again a smooth point blow up, and we prove 
in (\ref{defs.of.WB.claim}) that
every deformation of a weighted point blow up is 
again a weighted point blow up if $t\geq 3$.
Thus every deformation of a variety $X$ in  $\mathcal X$ is obtained by
deforming the points $p_0,\dots,p_n,q\in \p^n$
and also  the local coordinate system used for the 
weighted  blow up at $p_0$.
Since $p_0,\dots,p_n\in \p^n$ are in general position,
we can move their deformations  back to the
coordinate vertices by $\aut(\p^n)$;
hence we can  assume that the points
$p_0,\dots,p_n\in \p^n$ stay fixed in any deformation.

The point $q$ and the local coordinate system used for the 
weighted  blow up at $p_0$ however can deform nontrivially.
With these choices,  only the $(\c^*)^{s+t}$-action
remains of $\aut(\p^n)$.
\hfill\qed
\end{say}

\begin{say}[Choosing smooth moduli problems] \label{strong.exmps}
Using the explicit description of the
weighted blow ups given in (\ref{w.bu.defn})
we immediately obtain:

\begin{claim}\label{canonic.res} For every $X\in {\mathcal X}(s,t,d)$
the singular set $\sing X$ is isomorphic to
$\p^{s-1}$ and (Zariski locally) $X$ along $\sing X$ 
 is isomorphic to
\vglue12pt
\hfill $
\displaystyle{\a^{s-1}\times \a^{t+1}/\tfrac1{d}(1,(-1)^t).}$\hfill\qed
\end{claim}

These singularities are simple enough that one can write down
an explicit resolution for them, giving
``canonical'' resolutions $X^*\to X$ for every
$X\in {\mathcal X}(s,t,d)$. 
We do this in (\ref{quot.res.say}). 
Moreover, we prove that the local deformation theory 
of  $X^*$ is identical to the 
local deformation theory 
of  $X$.

By a suitable choice of the
polarization  $(X^*,H^*)$ we get a
smooth polarized moduli problem 
$({\mathcal X}^*(s,t,d), \mathcal H^*)$,
where the contraction $X^*\to X$ induces an
isomorphism of  the moduli spaces
${\mathcal X}^*(s,t,d)\cong {\mathcal X}(s,t,d)$.

\begin{prop} \label{versal.mod.exmps}
Notation and assumptions are as in
{\rm (\ref{torus.conds})} and {\rm (\ref{csmp.say}).}
 For $d\gg 1${\rm ,} there is an open subset
$$
{\mathcal X}^0(s,t,d)\subset {\mathcal X}(s,t,d)\cong {\mathcal X}^*(s,t,d)
$$
such that the $(\c^*)^{s+t}$-action is proper on
${\mathcal X}^0(s,t,d)$
and $U/T$ is isomorphic to a closed subscheme of
the quotient
$U/T\into {\mathcal X}^0(s,t,d)/(\c^*)^{s+t}$.

Thus ${\mathcal X}^0(s,t,d)/(\c^*)^{s+t}$ is a versal moduli
problem for smooth{\rm ,} polarized varieties as
in {\rm (\ref{vers.qp}.3)} which contains $U/T$ as a
closed subscheme.\hfill\qed
\end{prop}

All that remains is to find  examples
satisfying (\ref{torus.conds}) where $U/T$ is not quasi-projective.
\end{say}
 
\begin{exmp}\label{easy.exmps}
Consider $\a^{2t}$ with coodinates $y_1,\dots,y_t,x_1,\dots,x_t$.
Let $T$ be the torus $(\c^*)^t$ acting by
$$
y_i\mapsto \lambda_i\lambda^2_{i+1}y_i\qtq{(with $t+1= 1$), and}
x_i\mapsto \lambda_ix_i.
$$
Set $U_i:=(y_i\prod_{j\neq i}x_j\neq 0)$
and $U=\cup_iU_i$.
The $T$ action is free on $U$. \smallbreak %pagebreak

A polarization consists of an ample line bundle on
$\p^{2t}$, together with a linearization, that is, a
choice of the lifting of the $T$-action.
These correspond to characters
$$
\chi(b_1,\dots,b_t):(\lambda_1,\dots,\lambda_t)\mapsto 
\lambda_1^{b_1}\cdots \lambda_t^{b_t}.
$$
The $T$-equivariant monomials under this polarization are
of the form
$$
\bigl(\frac{y_1}{x_1x_{2}^2}\bigr)^{a_1}
\cdots
\bigl(\frac{y_t}{x_tx_{1}^2}\bigr)^{a_t}
\cdot 
\bigl(x_1^{b_1}\cdots x_t^{b_t}\bigr)^{m}.
$$
A $T$-orbit is semistable in the polarization given by 
$\chi(b_1,\dots,b_t)$ if and only if there is a monomial as above
which is  nonzero on the orbit.

Consider the orbit $C_i:=(x_i=0, y_j=0\ \forall\ j\neq i)$.
A monomial nonzero on $C_i$ can involve only
$y_i$ and the $x_j$ for $j\neq i$.
Thus, in the above form,  $a_j=0$ for $j\neq i$ and we have
a monomial of the form
$$
\bigl(\frac{y_i}{x_ix_{i+1}^2}\bigr)^{a_i}\cdot 
\bigl(x_1^{b_1}\cdots x_t^{b_t}\bigr)^{m}
$$
which does not contain $x_i$.
Thus  $a_i=mb_i$ and $2a_i\leq mb_{i+1}$, which
 is only possible if $b_{i+1}\geq 2b_i$.

Any collection of $t-1$ such inequalities has a common nonzero solution,
but all $t$ of them together lead to $b_1=\cdots=b_t=0$.

Thus we conclude that any $t-1$ orbits
in $U/T$ are contained in a quasi-projective open subset
but the $t$ orbits $C_1,\dots,C_t$ are
not contained in a quasi-projective open subset
of $U/T$. For $t\geq 3$ any 2 orbits are simultaneously stable
 with respect to some
polarization, so the quotient is separated.
Thus $U/T$ is  a variety which is not quasi-projective.
(For $t=2$ we get a nonseparated quotient.)
\end{exmp}

\begin{exmp}\label{simpelst.exmp}
The simplest proper but nonprojective toric variety $Y$
was found by Miyake and Oda, see \cite[\S 2.3]{oda}.
By \cite{cox}, this can also be obtained as the
quotient of an open subset of $\a^7$ by
a  $(\c^*)^4$-action.
I thank H.\ Thompson  for providing the following
explicit description.

Consider $\a^7$ with coodinates $y_1,y_2,y_3,x_1,x_2,x_3,x_4$.
Let $T$ be the 4-torus $(\c^*)^4$ acting by
\begin{multline*}
(y_1,y_2,y_3,x_1,x_2,x_3,x_4)\\ \mapsto
(\lambda_1\lambda_2\lambda_4y_1,\lambda_1\lambda_2\lambda_3y_2,
\lambda_1\lambda_3\lambda_4y_3,\lambda_1x_1,
\lambda_2x_2,\lambda_3x_3,\lambda_4x_4).
\end{multline*}
Let $U=\a^7\setminus Z$ where $Z$ is the subscheme corresponding to
the ideal
$$
(y_1,x_1)\cap(y_1,x_4)\cap(y_2,x_1)\cap(y_2,x_2)
\cap(y_3,x_1)\cap(y_3,x_3)\cap(x_2,x_3,x_4).
$$
Then $U/T$ is isomorphic to the  Miyake-Oda
example.

It is rather straightforward, though somewhat tedious, to check directly
that  $U/T$ is not projective
by looking at the set of stable points under all
possible polarizations.
\end{exmp}

\section{Weighted blow ups}

\begin{defn}\label{w.bu.defn}
 Let $x\in X$ be a smooth point on a variety of dimension $n$
and $(u_1,\dots,u_n)$ local coordinates. Choose positive integers
$(a_1,\dots,a_n)$, called weights.
This assigns weights to any monomial
by the rule
$$
w(u_1^{m_1}\cdots u_n^{m_n})=m_1a_1+\cdots +m_na_n.
$$
Let $I_c\subset \o_{x,X}$ be the ideal generated by
all monomials of weight at least $c$. We can also view $I_c$ as an ideal 
sheaf on $X$.
The scheme 
$$
B_{\mathbf u, \mathbf a}X:=
\proj_X\bigl(\textstyle{\sum_{c=0}^{\infty}}I_c\bigr)
$$
is called the {\it weighted blow up} of $X$
with coordinates $\mathbf u=(u_1,\dots,u_n)$ and weights
 $\mathbf a=(a_1,\dots,a_n)$.

In order to describe the local coordinate charts,
we use the notation 
$$
\a^n(u_1,\dots,u_n)/\tfrac1{d}(b_1,\dots,b_n)
$$
to denote the quotient 
of $\a^n$ with coordinates $u_1,\dots,u_n$
by the cyclic group of $d^{\rm th}$ roots of unity $\mu_d$  acting as
$$
\rho(\epsilon):(u_1,\dots,u_n)\mapsto 
(\epsilon^{b_1}u_1,\dots,\epsilon^{b_n}u_n).
$$
As a further shorthand, 
$$
\a^{s+t}/\tfrac1{d}(d^s,1^t)
$$
indicates that $s$ of the $b_i$ are $d$, and
$t$ of the $b_i$ are $1$.

With these conventions, (\'etale)
local coordinate charts on $B_{\mathbf u, \mathbf a}X$
are given by the quotients
$$
\a^n(x_{1,i},\dots,x_{n,i})/\tfrac1{a_i}(-a_1,\dots,-a_{i-1},1,
-a_{i+1},\dots,-a_n).
$$
The projection map is given by
\begin{gather*}
u_1=x_{1,i}x_{i,i}^{a_1},\dots, 
u_{i-1} = x_{i-1,i}x_{i,i}^{a_{i-1}},\\
u_i=x_{i,i}^{a_i},
u_{i+1}
 = x_{i+1,i}x_{i,i}^{a_{i+1}},\dots,
u_n=x_{n,i}x_{i,i}^{a_n}.
\end{gather*}

Let us now consider the special case when
$n=s+t$ and 
$(a_1,\dots,a_n)=(d^s,1^t)$.
Then the singular charts on $B_{\mathbf u, \mathbf a}X$
are of the form
$$
\a^{s+t}/\tfrac1{d}(1, (-1)^{t}, 0^{s-1}),
$$
proving (\ref{canonic.res}).

For weights $(d^s,1^t)$,
we get that
$$
I_c=(u_1,\dots,u_s)+m_x^c\qtq{for $c\leq d$,}
$$
and  the ideals $I_c$ are all determined by $I_d$. This in turn is
determined by the ideal $(u_1,\dots,u_s)$ modulo $m_x^d$.
Thus we conclude:

\begin{claim} The space $\mathcal W(s,t,d)$
of all weighted blow ups of weight
$(d^s,1^t)$ centered at a smooth point $x\in X$ 
can be identified with the subscheme of the
Hilbert scheme of points on $X$ parametrizing
the quotients $\o_{x,X}/I_d$.\hfill\qed
\end{claim}

Assume now that $X=\a^{s+t}$ with coordinates
$(y_1,\dots,y_s,x_1,\dots,x_t)$.

It is an open condition on $\mathcal W(s,t,d)$
 to assume
that the $y$-terms in the linear parts
of $u_1,\dots,u_s$ are linearly independent.
If this holds then by a linear change of coordinates we can
get a different generating set $u'_1,\dots,u'_s$ of the ideal
$(u_1,\dots,u_s)$ such that 
$$
u'_i=y_i+(\mbox{linear in $x$})+(\mbox{higher terms in $x,y$}).
$$
Then by a nonlinear coordinate change we can get
a final generating set $u^*_1,\dots,u^*_s$
such that 
$$
u^*_i=y_i+h_i(x_1,\dots,x_t)\qtq{where $h({\mathbf 0})=0$ and
$\deg h_i\leq d-1$.}
\leqno{(*)}
$$
The generators  $u^*_1,\dots,u^*_s$ are
uniquely determined
by the ideal $(u_1,\dots,u_s)$ modulo $m_{\mathbf 0}^d$
and hence by the weighted blow up.
So we conclude:

\begin{claim} There is an open subset 
 $\mathcal W^1(s,t,d)\subset \mathcal W(s,t,d)$
such that every weighted blow up in
$\mathcal W^1(s,t,d)$ can be uniquely given 
by local coordinates as in (*).\hfill\qed
\end{claim}
\end{defn}

The $(\c^*)^{s+t}$ action on $\mathcal W^1(s,t,d)$
is still pretty bad. To remedy this, we look at
a subset $\mathcal W^0(s,t,d)\subset \mathcal W^1(s,t,d)$
consisting of those local  coordinate systems
$u^*_i=y_i+h_i(x_1,\dots,x_t)$
such that 
$h_i(x_1,\dots,x_t)$ 
contains $\textstyle{\prod_j}x_j^{c_{ij}}$ with nonzero coefficient
for every $i$ where the $c_{ij}$ are as in (\ref{torus.conds}).

\begin{lem} With notation as above{\rm ,}  for $d>\max_i\{\sum_j c_{ij}\}${\rm ,}
\begin{enumerate}
\item $\mathcal W^0(s,t,d)$ is $(\c^*)^{s+t}$-invariant{\rm ,}
\item the diagonal $(\c^*)^{s+t}$-action  on
$\mathcal W^0(s,t,d)\times U$ is proper{\rm ,} and
\item the local coordinate system 
$u^*_i=y_i+\textstyle{\prod_j}x_j^{c_{ij}}$
is invariant under $T$ but not invariant under
any other element of  $(\c^*)^{s+t}$.
\end{enumerate}
\end{lem}

\Proof The action of $(\mu_1,\dots,\mu_s,\lambda_1,\cdots,\lambda_t)\in (\c^*)^{s+t}$
is given by
$$
y_i+h_i(x_1,\dots,x_t)
\mapsto \mu_iy_i+h_i(\lambda_1x_1,\dots,\lambda_tx_t)
\mapsto 
y_i+\mu_i^{-1}h_i(\lambda_1x_1,\dots,\lambda_tx_t).
$$
In particular, for $u^*_i=y_i+\textstyle{\prod_j}x_j^{c_{ij}}$
we get the action.
$$
y_i+\textstyle{\prod_j}x_j^{c_{ij}}\mapsto \mu_iy_i+
\left(\textstyle{\prod_j}\lambda_j^{c_{ij}}\right)
\left(\textstyle{\prod_j}x_j^{c_{ij}}\right)
\mapsto y_i+\left(\mu_i^{-1}\textstyle{\prod_j}\lambda_j^{c_{ij}}\right)
\left(\textstyle{\prod_j}x_j^{c_{ij}}\right).
$$
Thus we have invariance if and only if $\mu_i=\textstyle{\prod_j}\lambda_j^{c_{ij}}$
for every $i$, that is, only for elements of $T$.

Finally, the properness of the $(\c^*)^{s+t}$-action
is established in two steps.
First, note that $(\c^*)^{s+t}=T'+T$ where
$$
T':=\{(\mu_1,\dots,\mu_s,1,\dots,1)\in (\c^*)^{s+t}\}.
$$
Using the $T'$ action we can uniquely normalize
every coordinate system
$u^*_i=y_i+h_i(x_1,\dots,x_t)$ in $\mathcal W^0(s,t,d)$
to the form
$$
u^*_i=y_i+h_i(x_1,\dots,x_t)\ \mbox{where}\ 
\textstyle{\prod_j}x_j^{c_{ij}}\ 
\mbox{appears with coefficient 1.}
$$
Such coordinate systems form
a closed subset 
$\mathcal W^0_1(s,t,d)\subset \mathcal W^0(s,t,d)$,
and
$$
\bigl(\mathcal W^0(s,t,d)\times U\bigr)/(\c^*)^{s+t}\cong
\bigl(\mathcal W^0_1(s,t,d)\times U\bigr)/T.
$$
Since the $T$ action on $U$ is already proper
by assumption, it is also proper on
$\mathcal W^0_1(s,t,d)\times U$.
\hfill\qed

\section{Deformation and resolution of weighted blow ups}

\begin{say}[Deformation  of weighted projective spaces]
For an introduction to weighted projective spaces,
see \cite{dol}.

The infinitesimal
deformation space of a
weighted projective space  can be quite large.
The situation is, however, much better if
we assume that 
every singularity is a quotient singularity and
the singular set has codimension $\geq 3$.
In this case, by \cite{sch}, the local deformations
at every singular point are trivial, hence the first order
global deformations are  classified by 
$\ext^1(\Omega_X,\o_X)$.  On the
weighted projective space $X=\p(a_0,\dots,a_n)$
there is an exact sequence
$$
0\to \Omega_X\to \sum \o_X(-a_i)\to \o_X\to 0,
$$
and for $n\geq 3$ we obtain that $\ext^1(\Omega_X,\o_X)=0$.
Hence we conclude:

\begin{claim} If the singular set of $\p(a_0,\dots,a_n)$
has codimension at least 3, then 
every local deformation of $\p(a_0,\dots,a_n)$ is trivial.
In particular, every local deformation 
of $\p(1,1,1,a_3,\dots,a_n)$ is trivial.\hfill\qed
\end{claim}
\end{say}

\begin{say}[Deformation  of weighted blow ups]
\label{def.of.WB}

Let $X\subset Y$ be a proper subscheme of $Y$.
If  deformations of $X\subset Y$ are locally unobstructed, then
the obstructions to deforming $X$ in the Hilbert scheme 
$\hilb(Y)$ lie in  $$H^1(X, \Hom(I_X,\o_X)),
$$
where $I_X\subset \o_Y$ is the ideal sheaf of $X$.
If every singularity of $Y$ is a quotient singularity,
$X$ is a divisor  and
the singular set has codimension $\geq 3$,
then the deformations are locally unobstructed.

For the weighted blow up
$$
B_{(a_1,\dots,a_n)}X\to X\qtq{with exceptional divisor}
\p(a_1,\dots,a_n)\cong  E\subset B_{(a_1,\dots,a_n)}X
$$
the normal bundle is $\o(\sum a_i)$, hence
there are no obstructions for $n\geq 3$.
Thus $E$ deforms with any deformation of $B_{(a_1,\dots,a_n)}X$,
and the induced deformation of $E$ is trivial
if $a_1=a_2=a_3=1$. The exceptional divisor
$E$, its normal bundle and the local structure of
 $B_{(a_1,\dots,a_n)}X$
along the singular points determine 
$B_{(a_1,\dots,a_n)}X$ in a formal or analytic neighborhood
of $E$  (cf.\ \cite[Lem.~9]{hi-ro} or \cite[3.33]{mo-notnef}).
Thus every deformation of $B_{(a_1,\dots,a_n)}X$
is trivial in an analytic neighborhood of
the exceptional divisor. 
We conclude:

\begin{claim}\label{defs.of.WB.claim} If $a_1=a_2=a_3=1$ then every
deformation of $B_{(a_1,\dots,a_n)}\p^n$
is obtained by changing the local coordinate system
that defines the weighted blow up.\hfill\qed
\end{claim}

It is easy to express nontrivial deformations
of the weighted blow ups $B_{(1,d^{n-1})}\p^n$,
but rigidity probably holds assuming only that
$a_1=a_2=1$.
\end{say}

\begin{say}[Resolution of quotient singularities]
\label{quot.res.say}

It is not easy to get
resolutions for an  arbitrary cyclic quotient singularity, but for the
singularities\break  $\a^{t+1}/\frac1{d}(1,(-1)^t)$
there is a rather simple  resolution.

The key observation is that if we have the weighted blow up
$B_{(d-1,1^{t})}\a^{t+1}$ then the cyclic goup
action $\frac1{d}(1,(-1)^t)$ on $\a^{t+1}$
lifts to an action on\break $B_{(d-1,1^{t})}\a^{t+1}$
which acts trivially on the exceptional divisor.

Indeed, let us start with
$$
\a^{t+1}(x_0,\dots,x_n)/\tfrac1{d}(1,(-1)^t).
$$
The key chart  $U_0\subset B_{(d-1,1^{t})}\a^{t+1}(x_0,\dots,x_n)$ is
$$
U_0\cong \a^{t+1}(y_0,\dots,y_n)/\tfrac1{d-1}(1,(-1)^t),
\qtq{where $x_0=y_0^{d-1}$ and $x_i=y_iy_0$.}
$$
Thus we see that the $\frac1{d}(1,(-1)^t)$
action on
$\a^{t+1}(x_0,\dots,x_n)$ lifts to\break
$\a^{t+1}(y_0,\dots,y_n)$
as $\tfrac1{d}(-1,0^t)$.
Since
$$
\a^{t+1}(y_0,\dots,y_n)/\tfrac1{d}(-1,0^t)\cong 
\a^{t+1}(y_0^{d},y_1,\dots,y_n),
$$
we conclude that the quotient of $U_0$ is
$$
U_0/\mu_d\cong
\a^{t+1}(y_0^{d},y_2,\dots,y_n)/\tfrac1{d-1}(d,(-1)^t)=
\a^{t+1}/\tfrac1{d-1}(1,(-1)^t).
$$
Thus we get a partial resolution
$$
g_1:X_1\to X_0\cong \a^{t+1}/\tfrac1{d}(1,(-1)^t)
$$
whose exceptional divisor is $E_1\cong \p(d-1,1^{t})$
and $X_1$ has a unique singular point
of the form $\a^{t+1}/\frac1{d-1}(1,(-1)^t)$.
Moreover, the discrepancy (cf.\ \cite[\S 2.3]{km-book})
of $E_1$ is 
$a(E_1)=\frac{t-1}{d}$.

Working by induction, we thus obtain a tower
$$
X_{d-1}\stackrel{g_{d-1}}{\to} \cdots\stackrel{g_2}{\to}
X_1\stackrel{g_1}{\to} X_0\cong \a^{t+1}/\tfrac1{d}(1,(-1)^t)
$$
where $X_{d-1}$ is smooth and  $g_i:X_i\to X_{i-1}$
contracts a  single exceptional divisor
$E_i\cong \p(d-i,1^{t})$ to a point.
Moreover, we also obtain that $E_1$ has minimal discrepancy
among all exceptional divisors over 
the singular point of $ \a^{t+1}/\frac1{d}(1,(-1)^t)$.
(Indeed, any exceptional divisor other than $E_1$
is also exceptional over $X_1$; thus it either
lies over the
unique singular point and by induction has
discrepancy at least $\frac{t-1}{d-1}>\frac{t-1}{d}$,
or lies generically over the smooth locus of
$E_1$ and has discrepancy at least $1+\frac{t-1}{d}$.)
As in \cite[Prop.~6.5]{ko-nash2} this implies that
$g_1:X_1\to X_0$ is unique. Indeed, given any other
projective birational morphism $g'_1:X'_1\to X_0$
with a single exceptional divisor $E'_1$ of 
discrepancy $\frac{t-1}{d}$, 
the induced birational map $h:X_1\map X'_1$ is a
local isomorphism near the generic point of
$E_1$ since $E_1$ is the unique 
exceptional divisor of 
discrepancy $\frac{t-1}{d}$. Thus 
$h$ is an isomorphism by (\ref{mats-mum.lem})
since the\break $g_1$-ample divisor $-E_1$ is transformed into the
$g'_1$-ample divisor $-E'_1$.
This implies that in the setting of
(\ref{canonic.res}) the local resolution process automatically
globalizes.

(Note that 
this is a much stronger uniqueness than 
in the case of\break $B_{(a_1,\dots,a_n)}\a^n\to \a^n$, where
we have  uniqueness  only up to 
 a local coordinate change in
$\a^n$. For  $B_{(a_1,\dots,a_n)}\a^n\to \a^n$
the discrepancy is  $\sum a_i$. This is not minimal
unless all the $a_i=1$, and the usual blow up
is indeed unique.
It is a general rule that for  minimal discrepancy divisors 
we can expect  stronger
uniqueness results.)

As in (\ref{def.of.WB}) we also obtain that
every deformation of $X_i$ is trivial
for $t\geq 3$.

Putting all of these together, we conclude:

\begin{claim}\label{defs.res.comefrom.base} If $t\geq 3$ then every
deformation of the above constructed canonical resolution
of $B_{(d^s,1^t)}\p^{s+t}$
is obtained by changing the local coordinate system
that defines the weighted blow up.\hfill\qed
\end{claim}
\end{say}

\section{Open problems}

The above examples show that moduli spaces of smooth
polarized varieties can be complicated.
My guess is that in fact they have a universality property
with respect to subspaces.

\begin{conj} Let $G$ be a linear algebraic group acting properly
on a quasi-projective scheme $W$.
Then there are
\begin{enumerate}
\item  a projective space $\p${\rm ,} 
\item an open subset
$U\subset \hilb(\p)$ parametrizing smooth varieties
such that $\aut(\p)$ acts properly on $U${\rm ,}
\item a homomorphism $G\to \aut(\p)${\rm ,} and
\item a $G$-equivariant closed embedding $W\to U${\rm ,}
\end{enumerate}
such that the corresponding morphism
 $W/G\to U/\aut(\p)$ is a closed embedding.
\end{conj}

This naturally leads to the following question, which
is quite interesting in its own right.

\begin{ques} Which algebraic spaces can be written as
geometric quotients of quasi-projective schemes?
\end{ques}

The paper \cite{tot} contains a detailed review
and
a necessary and sufficient  condition in terms of
resolutions by locally free  sheaves. Nonetheless,
 the answer is not known even for 
normal schemes or 
smooth algebraic spaces.

\demo{Acknowledgments}  I thank N.\ Budur, D.\ Edidin, S.\ Keel, G.\ Schumacher, 
 H.\ Thompson
and B.\ Totaro
for useful comments, references  and suggestions.
Part of the work was done during the ARCC workshop
 ``Compact moduli spaces and birational geometry''.
Partial financial support was provided by  the NSF under grant numbers 
DMS02-00883  and DMS-0500198.

\end{document}